\input{psfig}
\documentstyle[12pt,diagram]{article}
\def\forces{\parallel\!\!\! -}

\def\hexnumber#1{\ifcase#1 0\or1\or2\or3\or4\or5\or6\or7\or8\or9\or
        A\or B\or C\or D\or E\or F\fi }

\font\teneuf=eufm10
\font\seveneuf=eufm7
\font\fiveeuf=eufm5
\newfam\euffam
\textfont\euffam=\teneuf
\scriptfont\euffam=\seveneuf
\scriptscriptfont\euffam=\fiveeuf
\def\frak{\fam\euffam \teneuf}

\font\tenmsx=msam10
\font\sevenmsx=msam7
\font\fivemsx=msam5
\font\tenmsy=msbm10
\font\sevenmsy=msbm7
\font\fivemsy=msbm5
\newfam\msxfam
\newfam\msyfam
\textfont\msxfam=\tenmsx  \scriptfont\msxfam=\sevenmsx
  \scriptscriptfont\msxfam=\fivemsx
\textfont\msyfam=\tenmsy  \scriptfont\msyfam=\sevenmsy
  \scriptscriptfont\msyfam=\fivemsy
\edef\msx{\hexnumber\msxfam}

\mathchardef\upharpoonright="0\msx16
\let\restriction=\upharpoonright
\def\Bbb#1{\tenmsy\fam\msyfam#1}

\def\qed{{\vcenter{\hrule height.4pt \hbox{\vrule width.4pt height5pt
 \kern5pt \vrule width.4pt} \hrule height.4pt}}}
\def\notin{{\in}\kern-5.5pt / \kern1pt}
\def\ok{\vbox{\hrule height 8pt width 8pt depth -7.4pt
    \hbox{\vrule width 0.6pt height 7.4pt \kern 7.4pt \vrule width 0.6pt height 7.4pt}
    \hrule height 0.6pt width 8pt}}
\def\nt{{\leq}\kern-1.5pt \vrule height 6.5pt width.8pt depth-0.5pt \kern 1pt}
\def\sd{{\times}\kern-2pt \vrule height 5pt width.6pt depth0pt \kern1pt}
\def\zp#1{{\hochss Y}\kern-3pt$_{#1}$\kern-1pt}

\def\LL{{\Bbb L}}

\def\PP{{\Bbb P}}
\def\QQ{{\Bbb Q}}

\def\TT{{\Bbb T}}
\def\11{{\Bbb 1}}

\font\capit=cmcsc10 scaled\magstep0

\newtheorem{theorem}{Theorem}[section]
\newtheorem{proposition}[theorem]{Proposition}
\newtheorem{corollary}[theorem]{Corollary}
\newtheorem{definition}{Definition}
\newtheorem{lemma}[theorem]{Lemma}

\newtheorem{question}{Question}
\newtheorem{fact}{Fact}

\title{Chains of End Elementary Extensions of Models of Set Theory}
\author{Andr\'es Villaveces}

\begin{document}

\maketitle

\begin{abstract}
\noindent
Large cardinals arising from the existence of arbitrarily long end elementary
extension chains over models of set theory are studied here. In particular,
we show that the large cardinals obtained that way (`Unfoldable
cardinals') behave as a `boundary'
between properties consistent with `V=L' and existence of indiscernibles.
We also provide an `embedding characterisation' of the unfoldable cardinals
and study their preservation and destruction by various different forcings.
\end{abstract}

\bigskip
\noindent
I wish to thank my advisor, Kenneth Kunen, for many helpful
conversations and comments, always full of interesting insights.
I also wish to thank Ali Enayat for very helpful discussions about some of the
topics treated in this article, as well as James Cummings for his insightful
comments, and Mirna D\v zamonja for carefully reading excerpts of this
article.

\setcounter{section}{-1}
\bigskip
\section{Introduction.}

\medskip
\noindent
Let $(M,E)$ be a model of set theory. A model of the same type $(N,F)$
{\bf end extends} $(M,E)$ iff for every $a\in M$, the sets
$a_E = \{ b\in M | b E a\}$ and $a_F = \{ b\in N | b F a\}$ are the same.
In other words, elements of $M$ are not enlarged by the extension
from $M$ to $N$. The structure $({\cal E}_M,\prec _e)$
consists of all the non-trivial end elementary extensions (for short,
{\bf eee}s) of $M$, ordered by the relation $\prec _e$. Properties related
to the non-triviality of $({\cal E}_M,\prec _e)$ (i.e., the {\bf existence}
of eees of $M$) and its connection to the
structure of $M$ were first studied by Keisler, Silver and Morley in
[KeMo 68] and [KeSi 70]. Then, among others, Enayat in [En 84], Kaufmann in
[Ka 83] and Villaveces in [Vi $\infty$] continued the study in that
direction.

\bigskip
\noindent
In this article, we concentrate on the study of the existence of various
types of {\bf chains} in ${\cal E}_{\frak A}$, where ${\frak A}$ is
an expansion of ${\cal R}(\kappa )$ by monadic predicates. 

\noindent
Our analysis of chains in ${\cal E}_{\frak A}$ is based on the study
of the following definition of a large cardinal property:

\medskip
\noindent
{\bf Definition~\ref{definition:unfold}:} $\kappa$ inaccessible is
{\bf unfoldable} iff $\forall
S\subset \kappa \forall \lambda \exists$ well-founded $N,\hat S$ such that
$({\cal R}(\kappa ),\in ,S)\prec _e (N,\in, \hat S)$ and $o(N)\geq \lambda$.

\medskip
\noindent
We prove that unfoldable cardinals are exactly those at which
${\cal E}_{({\cal R}(\kappa ),\in ,S)}$ has arbitrarily long chains,
for all $S\subset \kappa$ (Proposition~\ref{proposition:chains}).

\noindent
We prove that unfoldable cardinals are in between weakly compact and Ramsey
cardinals (Theorem~\ref{theorem:ramseysareunf}).
We compare the relative consistency strengths between unfoldables and
several cardinal inside those bounds, in contexts compatible with $V=L$ and
in the presence of large cardinals.

\noindent
Unfoldability relativises to $L$.
When $V=L$ holds, unfoldables are quite big cardinals: they are beyond
all the first $\Pi ^n_m$-indescribable cardinals, and (if they exist)
they are also beyond the first subtle cardinal. On the other hand, the
consistency of the existence of subtle cardinals is stronger than that
of unfoldable cardinals: we prove that if $\theta$ is subtle, then
${\cal R}(\theta )\models `\exists \kappa (\kappa$ unfoldable)'
(Theorem~\ref{theorem:subtle}).

\noindent
In the presence of Ramsey cardinals, the situation changes rather
dramatically: the strong indescribability properties of unfoldables
disappear (in that context, unfoldable cardinals are $\Pi ^1_2$-definable).
The first unfoldable is strictly in between the first weakly compact
and the first Ramsey cardinal. By contrast to the $V=L$ situation,
unfoldables are now rather small as large cardinals: In addition to
weak compactness, a partition property much weaker than Ramseyness already
implies the existence of unfoldable cardinals.
 
\noindent
We also provide the following characterisation of unfoldability
in terms of elementary embeddings.

\medskip
\noindent
{\bf Theorem~\ref{theorem:embed}.}
Let $\kappa$ be inaccessible, $\lambda \geq \kappa$. Then $\kappa$ is 
$\lambda$-unfoldable iff
$$\forall M \Bigl( M\mbox { transitive}, M \models ZF^-, |M| = \kappa \in M$$
$$\Rightarrow \exists j,N \left[ N\mbox { transitive}, \lambda \in N,
j:M\to N, crit(j)=\kappa, |j(\kappa)| = \lambda \right] \Bigr).$$

\medskip
\noindent
We then provide conditions on iterations that guarantee the
preservation of unfoldability. We also provide examples of non preservation
of unfoldability, and study the consistency strength of the failure of $GCH$
at an unfoldable cardinal. Using the embedding characterisation, we highlight
certain analogies between a natural strengthening of unfoldability
and strong cardinals.

\noindent
The following table summarises several results on the relative strength
of unfoldability.

\medskip

\begin{tabular}[h]{|c|c|c|}
\hline 
 $V=L$ & $0^\sharp$ exists & Ramseys  \\
       & no Ramseys    & exist  \\
\hline 
       &          & Ramsey \\
{\bf unfoldable} &      & $\uparrow$ \\
$\uparrow$   &      & $\kappa$ w. comp. $+\kappa \to
(\omega _1)^{<\omega _1}_{<\kappa}$ \\
$\pi ^n_m$ ($\forall n \forall m$) & {\bf unfoldable} & $\uparrow$ \\
$\uparrow$ &         & {\bf unfoldable} ($\Pi ^1_2$-def.) \\
           & $\uparrow$?  &  $\uparrow$ \\
weakly compact & weakly compact & weakly compact \\
\hline
\end{tabular}

\bigskip
\bigskip
\noindent
The arrows in this table work as follows:

\noindent
$x\to y$ means `the first $x$ cardinal is strictly less than the first $y$'.

\bigskip
\noindent
Some notational comments are important at this point: Following two different
traditions, we freely switch between the two notations
`$V_\kappa$' and `${\cal R}(\kappa )$' when we denote the set of objects of
the universe of rank less than $\kappa$. Given a model $M$, $o(M)$ denotes
the ordinal height of $M$.

\section{Unfoldable Cardinals.}
\subsection{Long Chains in ${\cal E}^{wf}_{({\cal R}(\kappa ),\in ,S)}$.}

\noindent
The following definition is central in this paper. It is a definition of
a large cardinal property that mentions objects of arbitrarily high rank
in the universe. It turns out to capture many connections between the
existence of arbitrarily long Chains in models of the form
${\cal E}^{wf}_{({\cal R}(\kappa ),\in ,S)}$, and other large cardinal
properties.

\begin{definition}\label{definition:unfold}
A cardinal $\kappa$ is {\bf unfoldable} if and only if it is inaccessible and
given any ordinal $\alpha$ and any $S \subset {\cal R}(\kappa )$,
there exists a well-founded $M\in {\cal E}_{({\cal R}(\kappa ),\in ,S)}$ such
that $o(M)\geq \alpha$.
\end{definition}

\bigskip
\noindent
This definition looks quite strong (it mentions objects of
arbitrarily high rank).
Nevertheless, unfoldable cardinals arise quite soon in the hierarchy of large
cardinals: all Ramsey cardinals are unfoldable, and depending on additional
axioms, unfoldables can be found below smaller cardinals, as will soon be
shown.

\medskip
\noindent
The structure of chains in $({\cal E}_{({\cal R}(\kappa ),\in ,S)},\prec_e)$
is extremely rich whenever $\kappa$ is an unfoldable cardinal: not only does it
contain objects of arbitrary height, but it also contains {\bf chains}
of arbitrary height (hence the name `unfoldable' for such $\kappa$). More
precisely, we have

\begin{proposition}\label{proposition:chains}
Let $\kappa$ be an unfoldable cardinal. Then, for any $S\subset
{\cal R}(\kappa )$,
for every ordinal $\alpha$,
$({\cal E}^{wf}_{({\cal R}(\kappa ),\in ,S)},\prec_e)$
contains chains of order type $\alpha$.
\end{proposition}

\noindent
{\bf Proof:} Fix $S\subset \kappa$. By the inaccessibility of $\kappa$, there
exists a club $C\subset \kappa$ such that
$$\gamma \in C \Longrightarrow ({\cal R}(\gamma),\in ,S\cap \gamma)\prec
({\cal R}(\kappa),\in ,S).$$
So, $({\cal R}(\gamma ),\in ,S\cap \gamma )_{\gamma \in C}$ is an eee-chain.
Let now $\beta >\kappa$. By the unfoldability of $\kappa$, there exist $M$
a transitive model, $\hat{S}$ and $\hat{C}$ such that
$$({\cal R}(\kappa ),\in ,S,C)\prec _e(M,\in ,\hat{S},\hat{C}),\mbox { and }
o(M) > \beta.$$
Clearly, $\kappa \in \hat{C}$. Then, $({\cal R}(\gamma )^M,\in ,\hat{S}\cap
\gamma)_{\gamma \in \hat{C}\setminus C}$ is a chain of order type $o(M)$ in
$({\cal E}^{wf}_{({\cal R}(\kappa ),\in ,S)},\prec_e)$. \hfill $\Box$

\begin{figure}[ht]\hspace*{35.00mm}
\psfig{figure=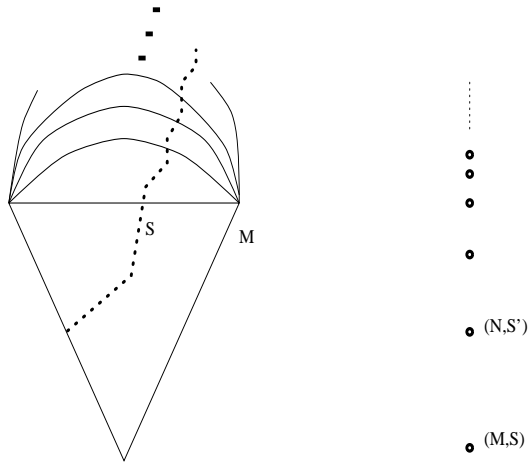,height=6cm,width=7cm}
\caption{ Unfoldables and Chains in
$({\cal E}_{({\cal R}(\kappa ),\in ,S)},\prec_e)$ } 
\end{figure}

\bigskip
\noindent
The parameter $S$ in the definition of unfoldable cardinals is crucial for
them to be of any interest: as Enayat points out in [En $\infty$], given any
model $M$ of $ZF$ such that ${\cal E}_M \not= 0$ (say $N\succ _eM$), there
exists (by a theorem of Schmerl) $N'\prec _eN$ such that cof$(N')=\omega$, and
$N'\in {\cal E}_M$. But then, by Keisler-Morley [KeMo 68], this $N'$ has
arbitrarily large eees. So, $M$ also has arbitrarily large eees. Thus, the
definition of unfoldable cardinals is nontrivial only when extending
models together with arbitrary additional
parameters is required. The consistency strength of the `lightface version'
of unfoldability is also quite low: if $V={\cal R}(\kappa )$, where $\kappa$
is the first Mahlo, then there are many inaccessible cardinals $\kappa$ in $V$
such that $({\cal E}_{{\cal R}(\kappa )},\prec_e)$ has an eee-chain of
length $ORD$. 

\medskip
\noindent
On the other hand, by iterating the $M$-ultrapower
construction over $M = {\cal R}(\kappa )$, for $\kappa$ weakly compact,
we obtain arbitrarily long chains of eees of $M$. They need not all be
well-founded: for example, if $V=L$, $\kappa ^+$ is an upper bound on the
length of chains of {\bf well-founded} eees of $M$. As Ali Enayat explained to
the author, a compactness argument, together with a forcing argument provides
$2^\kappa$ non-isomorphic eees of $M$, all of them well-founded, of the
same height, and of size $\kappa ^+$. These form a $2^\kappa$-antichain in
${\cal E}_M^{wf}$. The argument is a natural generalisation of the following
theorem, proved by Enayat in~[En $\infty$], and used by him
as part of the classification of countable models of complete extensions
of ZF. He proved that if $\tau$ is a countable linear order which is not
well-founded, and $T$ is a complete extension of ZF, if there exists a model
of $T$ with order type $\tau$, then there exist $2^{\aleph _0}$ of them.

\begin{theorem}(Enayat~[En $\infty$])
Let $M$ be a countable model of ZF. Then $M$ has continuum many non~isomorphic
eees.
\end{theorem}

\noindent
{\bf Sketch of the Proof:} first, one constructs a class $X$ in $M$ such
that $(M,X)\models ZFC(X)$ and such that every ordinal of $M$ is definable in $(M,X)$ (by forcing piecewise approximations to $X$ in
$M$). Then,
one constructs further expansions $(M,X,Y_\alpha )$, for $\alpha <2^\omega$,
such that $Y_\alpha \not= Y_\beta$ implies $(M,X,Y_\alpha ) \not \equiv
(M,X,Y_\beta )$. Finally, using Keisler-Morley ($M$ is countable!),
one gets eees $N_\alpha \succ _e M$, such that $N_\alpha$ codes
$(X,Y_\alpha )$, for $\alpha <2^\omega$. Then the $2^\omega$ different
$N_\alpha$s are pairwise non~isomorphic.
\hfill $\Box$

\bigskip
\noindent
The following summarises the previous discussion.

\begin{proposition} Let $\kappa$ be weakly compact. Then
\begin{description}
\item[ i ] ${\cal E}_{{\cal R}(\kappa )}$ has Ord-like chains,
\item[ ii ] ${\cal E}_{{\cal R}(\kappa )}^{wf}$ has chains of length
$\alpha$, for each $\alpha <\kappa ^+$,
\item[ iii ] $({\cal E}_{{\cal R}(\kappa )}^{wf})/_\approx$ has antichains
of size $2^\kappa$. \hfill $\Box$
\end{description}
\end{proposition}

\bigskip
\noindent
The process of iterating ultrapowers of models of
the form $({\cal R}(\kappa ),\in ,S)$, for $S\subset {\cal R}(\kappa )$,
whenever that is possible (for example, when $\kappa$ is measurable, or
by taking various sorts of ${\cal R}(\kappa )$-ultrapowers, when $\kappa$
is weakly compact, completely ineffable, etc.)
yields eee-chains of length $ORD$. This fact does not extend to every
unfoldable cardinal.

\begin{fact}
The unfoldability of $\kappa$ does not imply the existence of $ORD$-like chains
in the structure $({\cal E}^{wf}_{({\cal R}(\kappa ),\in ,S)},\prec_e)$.
\end{fact}

\noindent
We establish this fact by making the following observations. But we need
first the following definition.

\begin{definition}
An inaccessible cardinal $\kappa$ is {\bf long~unfoldable} if and only if
for every $S\subset {\cal R}(\kappa )$, the structure ${\cal E}^{wf}_{({\cal
R}(\kappa ),\in ,S)}$ has $ORD$-like chains.
\end{definition}

\noindent
In other words, besides requiring $({\cal R}(\kappa ),\in ,S)$ to have
arbitrarily long (chains of) eees, here we ask that it actually have
a long eee-chain. This requirement is not superfluous: for instance, if
$V=L$ holds, then the only proper class is $L$. Suppose that there is
an unfoldable cardinal, let $\gamma$ be the first one. If $\langle (L_{\gamma
_\alpha},\in ,S_\alpha)|\alpha \in ORD\rangle$ were an $ORD$-like chain in
$({\cal E}_{(L_\gamma ,\in ,S)},\prec _e)$, then its union ($=L$) would
be an elementary extension of $L_\gamma$. But this is impossible:
$L_\gamma$ believes that there are no unfoldables, and $L$ knows that
$\gamma$ is unfoldable. But $L\succ L_\gamma$!

\medskip
\noindent
On the other hand, in the presence of Ramsey cardinals, plus an
additional condition on the kinds of chains we have, this situation
changes strongly: if $\kappa$ is a Ramsey cardinal, and $\lambda <\kappa$
is unfoldable, then $\lambda$ is also long~unfoldable

\medskip
\noindent
To see this, we take any unfoldable cardinal $\lambda$, some $S\subset
\lambda$, and we long~unfold
the model $({\cal R}(\lambda ),\in ,S)$
in two steps: first, we consider some Ramsey cardinal $\kappa >\lambda$,
and we unfold $({\cal R}(\lambda ),\in ,S)$ {\bf up to} $\kappa$ in the
following strong sense:

\bigskip

$(*)^\kappa _\lambda \qquad \left \{
\begin{tabular}{lc}
for every
$S\subset \lambda$, there is a chain $\langle (M_\alpha ,\in , S_\alpha) |
\alpha < \kappa \rangle$ in\\
$({\cal E}^{wf}_{({\cal R}(\lambda ),\in ,S)},\prec _e)$, such that for every
$\alpha <\kappa$, $M_\alpha \in {\cal R}(\kappa )$.

\end{tabular}\right .$

\bigskip
\noindent
Our second step corresponds to applying stretching to that chain, which by
$(*)^\kappa _\lambda$ consists of objects
from ${\cal R}(\kappa )$. This way, we get at the end an $ORD$-like chain of
{\bf set} models in
$({\cal E}^{wf}_{({\cal R}(\lambda),\in ,S)},\prec _e)$.

\medskip
\noindent
We have just proved the

\begin{proposition}

If there are cofinally many Ramsey cardinals in the universe and for every
unfoldable $\lambda$ there exists a Ramsey $\kappa >\lambda$ such that
$(*)^\kappa _\lambda$,
\underline{then} all the unfoldable cardinals are long~unfoldable.
\hfill $\Box$

\end{proposition}

\medskip
\noindent
Thus, the consistency strength of `ZFC + unfoldability
and long~unfoldability are equivalent' is not greater than that of
the theory `ZFC + Ramsey
cardinals are cofinal in the universe + for every
unfoldable $\lambda$ there exists a Ramsey $\kappa >\lambda$ such that
$(*)^\kappa _\lambda$'.

\medskip
\noindent
{\bf Remark:} If $\theta$ is an ineffable cardinal, then there exists
$\kappa <\theta$, ${\cal R}(\theta )\models `\kappa$ is long~unfoldable'.
The proof of this fact is very similar to that of Theorem~\ref{theorem:subtle}
below, and is left to the reader. Thus, the consistency strength of the
existence of long~unfoldables is less than that of the existence of
ineffable cardinals,
and is thus still compatible with $V=L$. In a later section, we study in
greater detail the consistency strength of unfoldability under various
different extensions of ZFC.
The following chain of implications summarises the previous discussion:

$$\mbox {Con (ZFC + } \exists \mbox { ineffable)} \Longrightarrow
\mbox {Con (ZFC + } \exists \mbox { long unfoldable)} \Longrightarrow$$
$$\mbox {Con (ZFC + } \exists \mbox { unfoldable)}.$$

\bigskip
\noindent
Other facts worth noticing (parallel to those mentioned by Enayat
in [En $\infty$] regarding
{\it intermediate} extensions) about $({\cal E}_M,\prec_e)$
(for general $M$) are

\begin{description}
\item[  i  ] it is closed upwards under chains: given any chain, one can take
its union: it is clearly an element of $({\cal E}_M,\prec_e)$, by Tarski's
elementary chain theorem.
\item[  ii  ] it is also closed downward: let $\langle M_\alpha |\alpha
<\gamma \rangle$
be a sequence of eees of $M$ with $M_\beta \prec _eM_\alpha$ for $\alpha <
\beta <\gamma$ and let $M_\gamma = \bigcap _{\zeta \in \gamma}M_\zeta$.
Suppose that $M_\alpha \models \exists x\varphi (\vec{a},x)$, for some
formula $\varphi (\vec{y},x)$, and for a finite sequence $\vec{a}$ from
$M_\gamma$. Then, since every $M_\alpha$ contains the initial segment
$({\cal R}(\theta ))^{M_0}$, where $\langle M_0,\vec{a}\rangle \models
[\theta = \mbox { first }\zeta \mbox { such that }\exists x\in
{\cal R}(\zeta )\varphi
(\vec{a},x)]$, we have that for some $b\in M_\gamma$, $M_\alpha \models
\varphi (\vec{a},b)$. So, $M_\gamma \prec M_\alpha$, for each $\alpha
<\gamma$. Clearly, we also have $M_\gamma \succ _eM$; hence $M_\gamma$
is really a minimum for the chain $\langle M_\alpha |\alpha
<\gamma \rangle$ in $({\cal E}_M,\prec_e)$.

\end{description}

\noindent
It is reasonable to think that other structural features of the same
structure should be obtained using the unfoldability assumption, maybe
reinforced in some way.  

\bigskip
\noindent
We begin next the study of the consistency strength of
unfoldability. As the next few sections show, this consistency strength may
depend up to some point on the specific axioms of set theory considered.

\bigskip
\subsection{How big are the unfoldable cardinals?}
 
\begin{theorem}\label{theorem:ramseysareunf}
Let $\kappa$ be a Ramsey cardinal. Then $\kappa$ is unfoldable.
\end{theorem}

\noindent
{\bf Proof:} Let $S\subset {\cal R}(\kappa )$. Since $\kappa$ is Ramsey,
there are stationary many weakly compacts below it. Also,
${\cal R}(\kappa )$ supports theories of indiscernibles, on languages
$L_{\alpha ,\omega}$, for $\alpha <\kappa$. Our aim is to apply
stretching to get arbitrarily high eees of the {\bf whole} model
$({\cal R}(\kappa ),\in ,S)$. For this, it is enough to get a Silver theory of
indiscernibles for the expanded model $({\cal R}(\kappa),a)_{a\in
{\cal R}(\kappa )}$. This will guarantee endness of the elementary extensions.
Let us denote by $T\in
\Sigma \left[ {\frak A}\right]$ the fact that $T$ is a theory of
indiscernibles for
the structure ${\frak A}$, in the logic $L_{|{\frak A}|,\omega}$.
Now, in order to get the theory we need, we build a $\kappa$-tree $\TT$ of
theories of indiscernibles in the appropriate infinitary languages
for elementary submodels of $({\cal R}(\kappa ),\in ,S)$ by decreeing
$$(T_\beta)_{\beta <\alpha }\in \TT$$
\centerline {iff}
$$\forall \beta <\alpha
\left[ T_\beta \in \Sigma \left[ ({\frak A},\in ,S\cap|{\frak
A}|,a)_{a\in A}\right] ,{\frak A}\prec {\cal R}(\kappa ), |A|<\kappa \right].$$

\noindent
In this tree, there are branches of arbitrary lengths $<\kappa$, since
$\kappa$ is a limit of weakly compacts, and thus given any $\alpha <\kappa$
weakly compact, elementary submodels ${\frak A}$ of ${\cal R}(\kappa )$ of
height $\alpha$
have eees, which implies the existence of a branch of length at least $\alpha$
on top of the corresponding theory). But the ramification of $T$ is
small for
$\kappa$ ($2^{\lambda }<\kappa$, since $\kappa$ is strongly inaccessible).
So, $T$ is a $\kappa$-tree. Since $\kappa$ is weakly compact, $T$ must have
a path (of length $\kappa$): this path corresponds to a theory of
indiscernibles for the `union of the supports', $({\cal R}(\kappa),a)_{a\in
{\cal R}(\kappa )}$. So, we can apply stretching, and thus
get arbitrarily long chains of eees over $({\cal R}(\kappa ),\in ,S)$,
for any $S\subset \kappa$. \hfill $\Box _{\ref{theorem:ramseysareunf}}$

\medskip

\noindent
Also, by the Keisler-Morley characterisation of weakly compact cardinals,
we have immediately

\begin{proposition}
If $\kappa$ is unfoldable, then it is weakly compact.
\end{proposition}

\noindent
{\bf Proof:} It suffices to notice that unfoldability requires the existence
of arbitrarily high elements in $({\cal E}^{wf}_{({\cal R}(\kappa ),\in ,S)},
{\prec}_e)$,
whereas by Keisler-Morley [KeMo 68], $\kappa$ weakly compact amounts
to just having
$({\cal E}^{wf}_{({\cal R}(\kappa ),\in ,S)},{\prec}_e) \not=  0$.
\hfill $\Box$

\bigskip
\noindent
As we attempt to determine more precisely the consistency
strength of unfoldability, the answer begins to depend on the particular
model of set theory we work in. We will first study the case $V=L$.

\bigskip
\subsection{The Height Problem for Unfoldability.}

\bigskip
\noindent
In [Vi $\infty$], the relation between the height of a model $M$ and the
{\it existence \/} of end elementary extensions of $M$ was
studied. Several possible situations were analised there, and various
independence results were obtained. The main result there was that
it is consistent that
in most cases, inner models $M\subset {\cal R}(\kappa )$ of height $\kappa$
do not have end elementary extensions even if ${\cal R}(\kappa )$ has them.
More precisely, we proved the

\begin{theorem}\label{theorem:height}
 The theory `ZFC + $\exists \lambda G(\lambda )$ + $\forall \kappa (\kappa $
 inaccessible not weakly compact $\to$ NED($\kappa$))' is consistent relative
to the theory `ZFC + $\exists \lambda G(\lambda )$. $\Box$
\end{theorem}

\begin{figure}[h]\hspace*{50.00mm}
\caption{ The Existence Height Problem } 
\end{figure}

\noindent
Here, $NED(\kappa )$) means that there is a transitive $M_\kappa \models
ZFC$ such that $o(M)=\kappa$ and ${\cal E}_M=0$, and $G(\lambda )$ denotes
any large
cardinal property preserved under Easton-type extensions where the iteration
process is not carried too often (e.g. measurability, etc.).

\bigskip
\noindent
The corresponding `Height Problem for
Unfoldability' has a trivial answer:

\begin{fact}
If $M$ is a transitive model of $ZFC$, and
$o(M)=\kappa$, for $\kappa$ unfoldable, then given any $S\subset M$,
the structure $({\cal E}_{(M,\in ,S)},\prec _e)$ has elements of
arbitrarily high rank.
\end{fact}

\begin{figure}[ht]\hspace*{45.00mm}
\psfig{figure=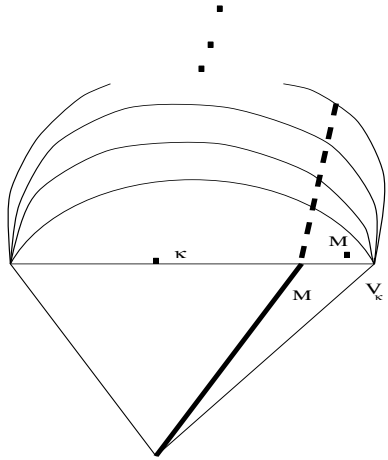,height=6cm,width=5cm}
\caption{ The Height Problem for unfoldables} 
\end{figure}

\noindent
To see this, just encode $M$ (as a subset of
${\cal R}(\kappa )$), and find (arbitrarily high ranked) eees of
the corresponding expanded structure. As before, this provides
the high ranked eees of $(M,\in ,S)$.

\noindent
This need not be true if the parameter is not required in the
definition: observe that the `weakened' form of
Unfoldability resulting from omitting the mention to the parameter $S$
in Definition \ref{definition:unfold} amounts to the requirement
${\cal E}_{\cal R(\kappa )}\not= 0$. In this case, the situation
is much more complicated; this was the main reason for the difficulty of
theorem \ref{theorem:height} in [Vi $\infty$].

\bigskip

\section{Unfoldables in L.}
\subsection{Relativisations.}

\medskip
\noindent
In this section, we study the relationship between Unfoldables in $V$ and
Unfoldables in $L$. We compare the two following situations: first,
properties such as relative consistency strength,
structure of $({\cal E}_{({\cal R}(\kappa ),\in ,S)},\prec_e)$, and
definability of Unfoldables, in the presence of axioms such as $V=L$, and
second, the same properties when certain large cardinals
are present.

\begin{theorem}
Let $\kappa$ be unfoldable in $V$. Then $\kappa$ is unfoldable in $L$.
\end{theorem}

\noindent
{\bf Proof:} Fix $S\subset {\cal R}(\kappa )$, $S\in L$. We need to build
$S'$ in $L$, in such a way that $[({\cal R}(\kappa ),\in ,S) \prec _e
(M,\in ,S')]^L$, for some $M$ high enough. Let $R \subset \kappa \times
\kappa$ code the construction of $S$ in $L$ (so that $\langle \kappa ,R\rangle
\approx \langle L(\xi ),\in \rangle$), for some $\xi$ such that $S\in
L(\xi )$, $\xi <\kappa ^+$. By the unfoldability of
$\kappa$ in $V$, we can get $\gamma$ of uncountable cofinality, and sets
$S',R'\subset L(\gamma )$ such that $(L(\kappa ),\in ,S,R)\prec _e(L(\gamma
),\in,S',R')$. $R'$ is well-founded, since $\gamma$ was taken of
uncountable cofinality. Then, it is possible to reconstruct $S'$ using
$R'$. Hence, $(L(\gamma ),\in ,S')$ is an eee of $(L(\kappa ),\in,S)$
in the sense of $L$. \hfill $\Box$

\medskip

\noindent
We remind to the reader the

\begin{definition}
An inaccessible cardinal $\kappa$ is {\bf subtle} if and only if given
any sequence $(S_{\alpha})_{\alpha <\kappa }$, where each $S_\alpha
\subset \alpha$, and given any $C$ club in
$\kappa$, there are $\beta <\gamma$, both in $C$, such that
$S_{\gamma }\cap \beta =S_{\beta }$.
\end{definition}

\noindent
This is a natural weakening of the notion of an ineffable cardinal.
Enayat [En 85] has obtained various results related to second-order
characterisations of end elementary extendability in terms of
completely ineffable cardinals.
Subtlety is a $\Pi _1^1$ property; yet, the first subtle is greater
than the first  $\Pi _n^1-indescribable$, for each $n$. They are thus
much larger than the first weakly compact.

\begin{theorem}\label{theorem:subtle}
Let $\theta$ be a subtle cardinal. Then, there exists $\kappa <\theta$
such that ($\kappa$ is unfoldable)$^{{\cal R}(\theta )}$.
\end{theorem}

\noindent
{\bf Proof:} Suppose by way of contradiction that the statement
$$\forall S\subset {\cal R}(\kappa )\forall \gamma <\theta \exists
\mbox { some } {\cal E}_{({\cal R}(\kappa ),\in 
,S)}\mbox {-chain of length }\gamma$$
fails for each $\kappa <\theta$. Then
we obtain the sequence
$\langle S_\kappa ,\gamma _\kappa ,\Delta _\kappa | \kappa <\theta \rangle$,
where for every $\kappa <\theta$, $\gamma _{\kappa }<\theta$, $S_\kappa
\subset {\cal R}(\kappa )$, $\Delta _\kappa$
is the elementary diagram of $({\cal R}(\kappa ),\in ,S_{\kappa})$,
and where all chains
of eees over $({\cal R}(\kappa ),\in ,S_{\kappa})$ are of length at most
$\gamma _\kappa <\theta$.

\noindent
Consider now a club $C$ such that if $\kappa \in C$,
then for all
$\alpha <\kappa$, $\gamma _{\alpha}<\kappa$. The subtlety of $\theta$
guarantees that we can pick $\kappa _1<\kappa _2$, both in $C$, such that
$$({\cal R}(\kappa _1),\in ,S_{\kappa _1})\prec ({\cal R}(\kappa _2),\in
,S_{\kappa _2}).$$
Now, this can be done in such a way that all the $\kappa$'s involved are
inaccessible. 
But then the set
$$C = \{ \alpha <\kappa _2 | ({\cal R}(\alpha ),\in
,S_{\kappa _2}\cap \alpha )\prec ({\cal R}(\kappa _2),\in ,S_{\kappa _2})\}$$
is a club in $\kappa _2$, and thus
$$\langle ({\cal R}(\alpha ),\in ,S_{\kappa _2}\cap \alpha)\rangle
_{\alpha \in C\setminus \kappa _1}$$
is an eee-chain over $({\cal R}(\kappa _1),\in ,S_{\kappa _1})$ of length
$\kappa _2$.
But this contradicts the fact that the maximum
length of such chains is $\kappa _1<\kappa _2$.
\hfill $\Box _{\ref{theorem:subtle}}$
\medskip

\begin{proposition}
If $\kappa$ is unfoldable, $\theta <\kappa$, and $(\theta \mbox{ is
unfoldable })^{\cal R(\kappa)}$, then $\theta$ is unfoldable in $V$.
\end{proposition}

\noindent
{\bf Proof:} Without loss of generality, take $\lambda >\kappa$, and let
$S\subset {\cal R}(\theta )$. Since ${\cal R}(\kappa ) \models
\forall \alpha \exists M_\alpha (o(M_\alpha ) \geq \alpha \wedge
{\cal R}(\theta )\prec _eM_\alpha )$, and $\kappa$ is at least
weakly compact, we can pick a sequence $C=\langle c_\alpha \rangle _{\alpha
<\kappa}$
of ordinals in $\kappa$ which encodes the `witnessing sequence'
$\langle M_\alpha \rangle _{\alpha <\kappa}$ in ${\cal R}(\kappa )$. Now,
we can pick a well founded eee of the structure $({\cal R}(\kappa ),\in
,{\cal R}(\theta ),S,C)$ of height $>\lambda$. Call this eee
$(N,\in ,{\cal R}(\theta )',S',C')$. Then $C'$ encodes a chain of
eees of ${\cal R}(\theta )' (={\cal R}(\theta ))$ of height
$o(N) > \lambda$ in $(N,\in ,{\cal R}(\theta ),S)$. Without loss of
generality, $N$ can have been picked high enough so that the last
fact holds as well in $V$. Since $\lambda$ was arbitrary, this provides
in $V$ eees of $({\cal R}(\theta ),\in ,S)$. \hfill $\Box$

\medskip

\noindent
So, the consistency strength of the existence of unfoldable cardinals lies
strictly between those
of the existence of weakly compact and of subtle cardinals. If $V=L$ holds,
then the first
weakly compact cardinal cannot be unfoldable: there are no chains of
well-founded models of length $\kappa ^+$ above ${\cal R}(\kappa )$
(otherwise, there would be such chains consisting of models of size $\kappa$;
but then ${\cal R}(\kappa ^+)$ would be an elementary extension of
${\cal R}(\kappa )$, which is clearly impossible). We also have that

\begin{proposition}\label{proposition:unsubtle}
If V=L holds, and there exist unfoldable and subtle cardinals, then the
first unfoldable is bigger than the first subtle.
\end{proposition}

\noindent
{\bf Proof:} Let $\kappa$ be unfoldable, and let $\alpha$ be subtle. Suppose,
without loss of generality, that $\alpha >\kappa$. Since $V=L$, there exists
$\delta >\alpha ^{++}$
such that $(L_\delta ,\in ,S')\succ _e(L_\kappa ,\in ,S)$. But then,
$L_\delta \models \exists \beta (\beta$ subtle), so $L_\kappa
\models \exists \beta (\beta$ subtle).
\hfill $\Box$

\medskip
\noindent
Contrast this to Theorem \ref{theorem:subtle}: below the
first subtle $\theta$, there are always unfoldables {\it in the sense
of\/ ${\cal R}(\theta )$}, yet are there NO real unfoldables below $\theta$.

\bigskip
\noindent
\section{Unfoldables, Ramseys and $\Pi _1^2$-sets.}

\noindent
Not only does the consistency strength of unfoldability vary radically when
we go from models of $V=L$ to models where larger cardinals are present:
the {\it definability\/} properties also change. We first look
at what happens in the extreme case $V=L$.

\begin{proposition} If V=L, and $\kappa$ is unfoldable, then for all m and n,
$\kappa$ is $\Pi _m^n$-indescribable.
\end{proposition}

\noindent
{\bf Proof:} Let $\varphi (v)$ be a $\Pi _m^n$ formula which holds `over'
${\cal R}(\kappa )$, with parameter $S$. Pick a high enough eee
of ${\cal R}(\kappa )$ (which necessarily is of the form $L(\lambda )$, by
well-foundedness and condensation). By reflecting it down,
we can get a club in $\kappa$ of inaccessibles at which $\varphi$ reflects
down. Note that parameters in the definition of unfoldables are not used in
this proof. \hfill $\Box$

\begin{corollary}\label{pi12def} If V=L holds, and there exist unfoldable and subtle
cardinals, then the
first unfoldable is bigger than $\pi _m^n$, for all $m$ and $n$.
\end{corollary}

\noindent
($\pi _m^n$ denotes the first $\Pi _m^n$-indescribable cardinal.)

\noindent
{\bf Proof:} Just look at the previous proposition and at Proposition
\ref{proposition:unsubtle}. \hfill $\Box$
\medskip

\noindent
{\bf Remark:} In the presence of Ramsey cardinals, the situation changes
radically:
unfoldable cardinals (below the first
Ramsey) are $\Pi _1^2$-definable via the formula

\bigskip
\qquad
$\varphi _{unf}(\kappa ) \equiv \left \{
\begin{tabular}{l}
$\forall S\subset \kappa \exists T[T$
is a theory of indiscernibles and \\
$(c_n)_{n<\omega}$ are indiscernibles and $T$ codes\\
the stretching of $Th(({\cal R}(\kappa ),\in ,S))$ up to $\kappa$ ].
\end{tabular}\right .$

\bigskip
\noindent
The formula $\varphi _{unf}$ is $\Pi _1^2$, and if $\kappa <$
first Ramsey, then $\varphi _{unf}(\kappa )$ holds if and only if $\kappa$
is unfoldable, by

\begin{lemma}
If $\theta$ is Ramsey and ${\cal E}_{({\cal R}(\kappa ),\in ,S)}$ has elements
arbitrarily high in ${\cal R}(\theta )$, then $\kappa$ is unfoldable.
\end{lemma}

\noindent
{\bf Proof:} applying Proposition \ref{proposition:chains} to
${\cal E}_{({\cal R}(\kappa ),\in ,S)}$ in ${\cal R}(\theta )$, get
cofinal chains of eees of $({\cal R}(\kappa ),\in ,S)$. Then stretch
them to arbitrary heights by using the fact that $\theta$ is Ramsey, and
a construction similar to the proof of Theorem \ref{theorem:ramseysareunf}.
\hfill $\Box$

\bigskip
\noindent
We are now in a position to compare the first unfoldable to the first Ramsey.
This is the main result in this section.

\begin{theorem}
Let $\kappa$ be the first unfoldable cardinal, and $\theta$ be the first
Ramsey. Then $\kappa <\theta$.
\end{theorem}

\noindent
But we actually prove the stronger
\begin{theorem}
Let $\kappa$ be weakly compact such that $\kappa \to (\omega _1)_{<\kappa }
^{<\omega }$. Then the first unfoldable is less than $\kappa$.
\end{theorem}

\noindent
(A Ramsey cardinal certainly satifies the two hypotheses of this
theorem; hence the previous Theorem.)

\noindent
{\bf Proof:} The proof will be divided into two steps. First, we will obtain
an indiscernible theory at $\kappa$. For reasons analogous to those given in 
Theorem~\ref{theorem:ramseysareunf}, getting this theory is enough to
make $\kappa$ unfoldable. In a second step, we show how to get unfoldables
below $\kappa$ from the result in the first step.

\medskip
\noindent
For each inaccesible cardinal $\alpha <\kappa$, we let $T_\alpha$ be a theory
in $\omega _1$ indiscernibles for the $L_{\alpha,\omega}$-theory of
$({\cal R}(\kappa ),\in ,a)_{a\in {\cal R}(\alpha )}$. These theories
are then pieced together by using the weak compactness of $\kappa$ in the
same way we did in the case of Theorem~\ref{theorem:ramseysareunf}.

\medskip
\noindent
We may now let $I$ be a set of indiscernibles for the model $\frak A =
({\cal R}(\kappa ),\in $, Skolem functions,$a)_{a\in {\cal R}(\kappa )}$,
with a {\bf least} $\omega$-th element among all the possible sets of
indiscernibles for that model. By the partition property $\kappa \to
(\omega _1)_{<\kappa }^{<\omega }$, we know that we can take $I$
of order type $\omega _1$. Let
now $M = \pi ({\frak H}_{\frak A}(I))$
(the Mostowski transitive collapse of the Skolem hull of $I$ in $\frak A$),
and let $I^*$ be the set of indiscernibles for
$M$ corresponding to $I$. Then
$o(I^*)=\omega _1$, and $I^*$ is a club in (the true) $\omega _1$. Also,
$$\alpha \in I^* \mbox { implies }[\alpha \mbox {
is unfoldable }]^M:$$
if not, then none of the elements of $I^*$ would be unfoldable according to
$M$. So, there would be counterexamples of the form
$S_\alpha \subset \alpha$ for each $\alpha \in I^*$. But by Fodor's
Lemma, there exists a club $C$ of $I^*$ such that all the
counterexamples corresponding to indiscernibles from $C$ necessarily
cohere (by the original choice of $I$
as the set of indiscernibles with the least possible $\omega$-th element
among them, and a remarkability argument).
But then, by elementarity, $$\alpha \in I \mbox { implies }[\alpha \mbox {
is unfoldable }]^{{\cal R}(\kappa )}.$$

\medskip
\noindent
Up to now, we have cardinals $\alpha$ that unfold up to $\kappa$.
By stretching the indiscernibles $I$ as in previous constructions,
we observe that $\alpha$ is also unfoldable in $V$. \hfill $\Box$

\bigskip
\noindent
The previous results give `upper bounds' for consistency of Unfoldability.
We now concentrate on the `lower bounds'. We already know that, under
$V=L$, the first unfoldable is at least the first completely
indescribable cardinal (if they exist), and is thence much greater than the
first weakly compact cardinal.

\medskip
\noindent
In general, as observed before, one cannot guarantee that the `gap' between
weakly compacts and unfoldables is as big as under $V=L$: in the
presence of Ramsey cardinals, unfoldables are $\Pi ^2_1$-definable. Still,
the fact that the first unfoldable is greater than the first weakly
compact may hold under some additional large cardinal hypotheses. One
important case is given by the following theorem.

\begin{theorem}\label{theorem:wcnotunf}
If there exists a Ramsey cardinal, then the first
unfoldable cardinal is greater than the first weakly compact. (And of course,
these two exist!)
\end{theorem}

\noindent
{\bf Proof:} Let $\lambda$ be Ramsey, and let $\kappa <\lambda$ be unfoldable (the
previous results show that this situation is possible). We prove that
$\kappa$ is not {\it the first\/} weakly compact cardinal. To reach a
contradiction, assume $\kappa =$ first weakly compact. For such a $\kappa$,
there is a chain in $({\cal E}_{({\cal R}(\kappa),\in ,S)},\prec _e)$ of
length $\kappa$ that cannot be extended to a chain of length $\kappa ^+$ in
the same structure. (It could still have {\it other\/} chains of arbitrarily
high length; we have not yet reached a contradiction). Now, let
$\vec M=\langle (M_\alpha ,\in ,S_\alpha )|\alpha <\theta \rangle$ be one
such chain, where $\theta$ is some ordinal less than $\kappa ^+$ and
$|M_\alpha |=\kappa$, for all $\alpha <\theta$. The chain
$\vec M$, being of length $<\kappa ^+$, can be coded by some
$Q\subset \kappa$. By the
unfoldability of $\kappa$, the model $({\cal R}(\kappa ),\in ,S,Q)$ has
some transitive end elementary extension $(N,\in ,\tilde S,\tilde Q)$
with $o(N)\geq \lambda$. But then, $\tilde Q$ codes an eee chain of 
length $o(N)$. Stretching the model $(N,\in ,\tilde S,\tilde Q)$: this
provides eee-chains of arbitrary length over $({\cal R}(\kappa),\in
,S)$. \hfill $\Box _{\ref{theorem:wcnotunf}}$

\bigskip
\noindent
The following table summarises several of the results of the last sections.

\bigskip
\begin{tabular}[h]{|c|c|c|}
\hline 
 $V=L$ & $0^\sharp$ exists & Ramseys  \\
       & no Ramseys    & exist  \\
\hline 
       &          & Ramsey \\
{\bf unfoldable} &      & $\uparrow$ \\
$\uparrow$   &      & $\kappa$ w. comp. $+\kappa \to
(\omega _1)^{<\omega _1}_{<\kappa}$ \\
$\pi ^n_m$ ($\forall n \forall m$) & {\bf unfoldable} & $\uparrow$ \\
$\uparrow$ &         & {\bf unfoldable} ($\Pi ^1_2$-def.) \\
           & $\uparrow$?  &  $\uparrow$ \\
weakly compact & weakly compact & weakly compact \\
\hline
\end{tabular}

\bigskip
\bigskip
\noindent
The arrows in this table work as follows:

\noindent
$x\to y$ means `the first $x$ cardinal is strictly less than the first $y$'.

\noindent
This table also reflects how the consistency strength of unfoldability
`decreases' when one moves from $V=L$ to the presence of larger cardinals.
The question mark corresponds to the fact we do not know yet if
the strict inequality between the first weakly compact and the first
unfoldable that holds in the two `extreme' cases also holds in the
intermediate column.

\bigskip
\noindent
\section{Preservation of Unfoldability.}

\noindent
Our main tool in the study of the relationship between iterated forcing
and unfoldability is the embedding characterisation from
Section~\ref{embchar}. We address the general problem of constructing
iterations that preserve unfoldability, and then study the specific case
of adding or destroying GCH. Several connections between unfoldable
cardinals and strong cardinals are encountered on the way.

\subsection{An embedding characterisation.}\label{embchar}

\medskip
\noindent
In this section, we provide a characterisation
of unfoldability in terms of elementary embeddings. This is a useful
version of the definition of these cardinals, and provides interesting
variants. On the other hand, our characterisation highlights a sharp
analogy between unfoldable and strong cardinals. This analogy suggests
natural paths for further research. We intend to use master condition
arguments in order to prove the preservation
of the unfoldability of cardinals under various iterations. We will
provide conditions under which elementary embeddings
$j:M\to N$, where $M$ and $N$ model
fragments of ZFC, lift to elementary embeddings $j:M[G]\to N[H]$, where $G$ is
$\PP$-generic over $M$, for some $\PP \in M$, and $H$ is $j(\PP )$-generic
over $N$. In what follows, we will refer to elementary
embeddings just as `embeddings', and in most cases, will abuse notation
by using `$j$' or `$i$' for both the basic embedding and its lifting.
Remember that a sufficient
condition for the lifting to exist is that
$$ \forall p \in G \mbox {  } j(p)\in H.$$

\noindent
Here is a reformulation of unfoldability in terms of elementary embeddings.
James Cummings suggested this line of research.

\begin{theorem}\label{theorem:embed}
Let $\kappa$ be inaccessible, $\lambda \geq \kappa$. Then $\kappa$ is 
$\lambda$-unfoldable iff
$$\forall M \Bigl( M\mbox { transitive}, M \models ZF^-, |M| = \kappa \in M$$
$$\Rightarrow \exists j,N \left[ N\mbox { transitive}, \lambda \in N,
j:M\to N, crit(j)=\kappa, j(\kappa) \geq \lambda \right] \Bigr).$$
\end{theorem}

\begin{figure}[ht]\hspace*{35.00mm}
\psfig{figure=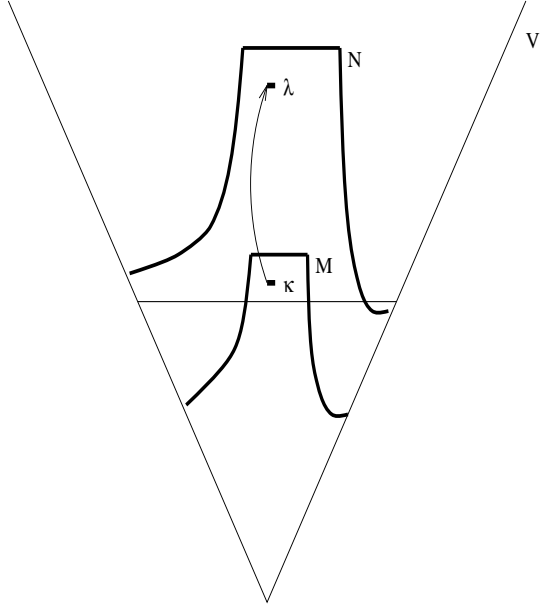,height=8cm,width=7cm}
\caption{ The embedding characterisation} 
\end{figure}

\noindent
{\bf Proof:} We first assume that $\kappa$ is $\lambda$-unfoldable.
Start with $M$ a transitive model of a finite fragment of ZFC,
such that $|M|=\kappa \in M$ and $M^{<\kappa}\subset M$. Code $(M,\in )$
by $(\kappa ,E)$ in such a way that the Mostowski collapse $\pi$ of
$(\kappa ,E)$ is such that
$\pi ``\kappa =M$, $\pi (0)=\kappa$. Then $E\subset \kappa ^2\subset V_\kappa$.
By the $\lambda$-unfoldability of $\kappa$, get $({\frak A},\in ,\tilde
E)\succ _e({\cal R}(\kappa),\in ,E)$, with $o({\frak A})\geq \lambda$. Now
interpret $\tilde E$ in $\frak A$, and take its transitive collapse $\tilde M$,
via $\tilde \pi$, so that $\tilde \pi (0) = \lambda$. Let also $j =\tilde
\pi \circ \pi^{-1}:M\to \tilde M$. Then $j(\kappa ) =
\pi \circ \pi^{-1}(\kappa ) = \tilde \pi (0) = \lambda$ and $|N| =
|\tilde M| = \lambda$, and we get the next commuting diagram.

\begin{figure}[ht]
$$
\setlength{\dgARROWLENGTH}{5.0em}
\begin{diagram} [\strut A]
\node{M} \arrow{e,t}{j} \node{N=\tilde M} \\
\node{(\kappa ,E)} \arrow{n,l}{\pi}\arrow{e,t}{\prec}\node{(\lambda,\tilde E)}
\arrow{n,r}{\tilde \pi}
\end{diagram}
$$
\end{figure}

\noindent
Let now $j =\tilde \pi \circ \pi^{-1}:M\to \tilde M$. Then $j(\kappa ) =
\pi \circ \pi^{-1}(\kappa ) = \tilde \pi (0) = \lambda$ and $|N| =
|\tilde M| = \lambda$.

\medskip
\noindent
For the converse, let $S\subset {\cal R}(\kappa )$, and take a model
of enough set theory $M$ such that $|M| = \kappa$, $M\supset
{\cal R}(\kappa )$, and both $\kappa$ and $S$ are elements of $M$. Obtain
$j$ and $N$ as in the statement of the theorem. Define
$$({\frak A},\in ,\tilde S) = (N_{j(\kappa )},\in ,j(S)).$$
Then
\begin{description}
\item[ i ] $({\cal R}(\kappa),\in ,S)\subset _e({\frak A},\in ,\tilde S)$:
membership in $N$ is just $\in$, and $\tilde S = j(S) = S\cup T$, for some
$T$ {\bf above\/} ${\cal R}(\kappa)$
(that is, $j(S)\cap {\cal R}(\kappa) = S$).
\item[ ii ] $({\cal R}(\kappa),\in ,S)\prec ({\frak A},\in ,\tilde S)$:
if $({\frak A},\in ,\tilde S)\models \exists x\varphi(x,a_1,\dots ,a_n)$, for
$a_i\in {\cal R}(\kappa)$, then
$$N\models \left[({\cal R}(j\kappa ),\in
,j(S))\models \exists x\varphi(x,j(a_1),\dots ,j(a_n))\right],$$
since $j(a_i) = a_i$ for each $i$. By elementarity, we have
$$M\models \left[({\cal R}(\kappa ),\in
,S)\models \exists x\varphi(x,a_1,\dots ,a_n)\right].$$
Finally, since $M_\kappa = {\cal R}(\kappa )$,
$$({\cal R}(\kappa ),\in
,S)\models \exists x\varphi(x,a_1,\dots ,a_n).$$
\item[ iii ] $o({\frak A}) \geq \lambda$, as $|j(\kappa )| = \lambda$
and $o(N_{j(\kappa )}) = j(\kappa )$.
\end{description}
\hfill $\Box_{\ref{theorem:embed}}$

\bigskip
\noindent
These results are a generalisation of Keisler-Silver's
characterisation of weak compactness as the existence of non-trivial
eees for all $S\subset \kappa$, and Kunen's characterisation of
weak compactness of $\kappa$ via the existence of an $M$-ultrafilter
on $\kappa$. (All these, for $\kappa$ inaccessible!)
A similar generalisation, but in the direction of $\Pi^m_n$-indescribability
was studied by Hauser in~[Ha 91]. We quote his theorem next
to observe the analogy, and also because we use it later.

\begin{theorem}(Hauser~[Ha 91])\label{hauser}
Let $m\geq 1, n\geq 1$. An inaccessible cardinal $\kappa$ is
$\Pi^m_n$-indescribable iff
$$\forall M\Bigl( M\mbox { trans} \wedge M\models ZF^- \wedge |M|=\kappa \wedge
\kappa \in M \wedge M^{<\kappa}\subset M$$
$$\Rightarrow \exists j,N \Big[ N\mbox { trans}\wedge |N|=|V_{\kappa+m-1}|
\wedge N \mbox { is }\Sigma ^m_{n+1}\mbox { -correct for }\kappa$$
$$\wedge \mbox { } j:M\to N \wedge crit(j)=\kappa \Big] \Bigr).$$

\end{theorem}

\bigskip
\noindent
{\bf Remarks: 1)} Our characterisation of $\lambda$-unfoldability
provides an analogy between these and $\lambda$-strong cardinals,
parallel to the analogy observed
by Kai Hauser in [Ha 91, p.445] for $\Pi^1_m$-indescribable cardinals and
$m$-hypermeasurable cardinals. In our case,
the $\lambda$-unfoldability of a cardinal can be
construed as an analogue of $\lambda$-strong cardinals in the same way as
weakly compact cardinals relate to measurable cardinals. In the following
table, the columns are the related pairs smaller/larger cardinal.

\bigskip
\qquad \qquad \begin{tabular}{|c|c|c|}
\hline
$\Pi ^1_m$-indescr.&weakly compact&$\lambda$-unfoldable\\
\hline
$m$-hypermeas.&measurable&$\lambda$-strong\\
\hline
\end{tabular}

\medskip
\noindent
{\bf 2)} Taking $\lambda$ arbitrarily large, the previous analogy extends to
one between unfoldable and strong cardinals.

\medskip
\noindent
{\bf 3)} The proof of Theorem \ref{theorem:embed} suggests a
strengthening of the
definition of $\lambda$-unfoldable cardinals: in the characterisation via
embeddings, we could require that $N$ contain the whole ${\cal R}(\lambda )$.
We get the following definition.

\begin{definition}
$\kappa$ inaccessible is $\lambda${\bf -strongly-unfoldable} if
and only if 
$$\forall M \Bigl( M\mbox { transitive}, M \models ZF^-, |M| = \kappa \in M,
M^{<\kappa}\subset M$$
$$\Rightarrow \exists j,N \Big[ N\mbox { trans},
\underbrace{{\cal R}(\lambda )\cup
\{ \lambda \} \subset N}_{\mbox {the new ingredient}},
j:M\to N, crit(j)=\kappa, j(\kappa) \geq \lambda \Big] \Bigr).$$
\end{definition}

\medskip
\noindent
It turns out that $\lambda$-strongly-unfoldable cardinals are $\Pi
^n_m$-indescribable for all $n$ and $m$: since $N \supset V_\lambda$,
it is $\Sigma ^n_m$-correct for all $n$ and $m$, and then Hauser's
characterisation of $\Pi ^n_m$-indescribability works here.

\medskip
\noindent
If $V=L$, then strong unfoldability and unfoldability are just the same.
On the other hand, strong unfoldability must be stronger than unfoldability
in the
presence of Ramsey cardinals, since in that case, as observed right after
Corollary \ref{pi12def}, the formula

\bigskip
\qquad
$\varphi _{unf}(\kappa ) \equiv \left \{
\begin{tabular}{l}
$\forall S\subset \kappa \exists T[T$
is a theory of indiscernibles and \\
$(c_n)_{n<\omega}$ are indiscernibles and $T$ codes\\
the stretching of $Th(({\cal R}(\kappa ),\in ,S))$ up to $\kappa$ ].
\end{tabular}\right .$

\bigskip
\noindent
{\bf 4)} This embedding characterisation of $\lambda$-unfoldability also
provides an extender characterisation, which may be useful if
finer approximations of the embedding in question are to be used: an
inaccessible cardinal $\kappa$ is $\lambda$-unfoldable if and only if
for all transitive $M\models ZF^-$, of size $\kappa$, closed under
$(<\kappa )$-sequences, with $\kappa \in M$, there exists an $M$-extender
(the natural relativisation of the notion of an extender to $M$) with
support $V_\lambda$.

\subsection{Iterations that preserve unfoldability.}
\noindent
To preserve unfoldability under iterations is not a trivial matter. Several
of the iterations which preserve other large cardinal properties such as
measurability, $\Pi ^1_n$-indescribability, ineffability, fail to
preserve unfoldability. Why? Let us first state a theorem
due to Levinski, and then explain why unfoldability fails under iterations
that work in other contexts.

\begin{theorem}(Levinski~[Le 95])\label{levinski}
Let $V\models ZFC$, and let $\PP _\kappa$ be a Reverse Easton iteration
of L\'evy collapses $\LL ((a_\alpha)^+,2^{a_\alpha})$, where $a_{\alpha +1}$
is forced over $\PP_{\alpha + 1}$ to be $(a_\alpha)^+$, and for limit
$\alpha$, $a_\alpha$ is forced over $\PP _\alpha$ to be the supremum of
all the $a_\xi$, for $\xi <\alpha$.

\noindent
\underline {Then}, if $V\models `\kappa$ is measurable (or
$\Pi ^1_n$-indescribable, or ineffable, or completely ineffable)', then
$V^{\PP _\kappa}\models `\kappa$ is measurable (or
$\Pi ^1_n$-indescribable, or ineffable, or completely ineffable)'.
\end{theorem}

\noindent
One would expect an analogous result for unfoldability. In the case one
does a Reverse Easton iteration of the usual forcings that `add new Cohen
subsets' to all inaccessibles
below an unfoldable $\kappa$, it may well happen that $\kappa$'s
unfoldability is $\bf not$ preserved. This is the content of the following
counterexample to the preservation of unfoldability.


\noindent
\begin{theorem}
Suppose that $V=L$. Let $\kappa$ be the {\bf first}
unfoldable cardinal. Then, the Reverse Easton Iteration described next
does not preserve $\kappa$'s unfoldability.
\end{theorem}

\noindent
{\bf Proof:} Since $V=L$ and $\kappa$ is the first unfoldable,
$({\cal R}(\kappa ),\in ,S)$ must
unfold to well-founded models of the form $({\cal R}(\beta ),\in ,\hat S)$.

\medskip
\noindent
{\bf Claim:} The previous $\beta$ cannot be an inaccessible cardinal.
Actually, even more is true: $\beta$
cannot be a {\bf cardinal}. The reasons are that $\kappa$ is the first
unfoldable, and the elementarity of the extension.

\medskip
\noindent
Assume furthermore that $\lambda$ is larger than $\kappa$, and
Easton-iterate the forcings that add new Cohen
subsets to every inaccessible $\alpha <\lambda$, to get the forcing
notion $\PP$. Assume also that there is some
inaccessible $\theta >\lambda$. Let $G$ be $\PP$-generic over $V$.

\noindent
Then, in $V[G]$, if $[{\cal R}(\kappa )]^{V[G]}$ is unfolded to some
model $M$ of
height $> 2^\theta$, $M$ will `think' that a generic subset of $\theta$ has
been added, since $\theta$ is inaccessible in $M$. But the forcing $\PP$
is a $(\lambda +1)$-iteration, and thus could not have added new unbounded
subsets to $\theta >\lambda$!
\hfill $\Box$

\bigskip
\noindent
In the next section, we will encounter a situation similar to the one
described in the last theorem, 

\bigskip
\noindent
On the other hand, there is a positive situation for preservation of
unfoldability under iterations. We could roughly describe the next theorem
as asserting that `$(<\kappa )$-support iterations of forcings of size
$<\kappa$ that satisfy additional lifting conditions preserve the
unfoldability of cardinals $\kappa$'. Although the `lifting conditions' will be
lengthy to describe, they turn out to be quite natural, as the reader may
check.

\begin{theorem}\label{presunf}
Suppose that $\kappa$ is $\lambda$-unfoldable, there are inaccessibles
above $\kappa$, and $\PP = \PP _{\kappa
+1}$ is an iteration of length $\kappa +1$. Suppose furthermore that
$\PP$ has $(<\kappa )$-supports, and that $|\PP _\kappa|=\kappa$, and
that $\PP$ satisfies the lifting conditions described in what follows next.
\underline{Then}, $V^\PP \models `\kappa \mbox { is }\lambda$-unfoldable'.
\end{theorem}

\bigskip
\noindent
We begin by describing our additional conditions.

\medskip
\noindent
{\bf The Lifting Conditions for Theorem \ref{presunf}:}
\begin{itemize}
\item[ \bf i ] $p\forces _\PP `\dot S \subset {\cal R}(\kappa )$', $\gamma$
is greater than the least inaccessible above $\kappa$ and such that
cof($\gamma )>\kappa$ and $V_\gamma \models ZF^-$.
\item[ \bf ii ] $\PP$ can be defined in ${\cal R}$(least inacc $>\kappa$),
so that $\PP ^{V_\gamma} = \PP$.
\item[ \bf iii ] (For $M\models ZF^-$ transitive
of size $\kappa$ such that $M^{<\kappa}\subset M$, $\kappa \in M$, and
for $i:M\to V_\gamma$ with critical point $(\kappa ^+)^M$.) Let
$\dot S^M$ be such that $i(\dot S^M)=\dot S$. Notice that $\dot S^M$ is a
$\PP ^M$-name. We define $i$ as a transitive collapse of $M$ with
ran($i)\subset V_\gamma$. This is possible by Mostowski's theorem in
$V_\gamma$. Wlog, $\dot S\in$ ran($i$). As $M^{<\kappa}\subset M$, the critical
point of $i$ is $(\kappa^+)^M$.
\item[ \bf iv ] (By the $\lambda$-unfoldability of $\kappa$, choose $j$ with
$|N|=\lambda \in N$, $N$ $\Sigma ^1_0$-correct for $\kappa$, and
$j:M\to N$, crit($j$)=$\kappa$. Wlog, $N$ can be picked in $V_\gamma$.)
There exist $G\ni p$, $\PP$-generic over $V_\gamma$,
$H$ $\PP
^M$-generic over $M$ and $K$ $(\PP _{j(\kappa )+1})^N$-generic over $N$
such that both $i$ and $j$ lift through the corresponding forcings

\begin{figure}[ht]
$$
\setlength{\dgARROWLENGTH}{5.0em}
\begin{diagram} [\strut A]
\node{M[H]} \arrow{e,t}{i} \node{V_\gamma [G]}
\node{M[H]} \arrow{e,t}{j} \node{N[K]}\\
\node{M} \arrow{n,l}{\PP ^M_{\kappa +1}}\arrow{e,t}{i}\node{V_\gamma}
\arrow{n,r}{\PP _{\kappa +1}}
\node{M} \arrow{n,l}{\PP ^M_{\kappa +1}}\arrow{e,t}{j}\node{N,}
\arrow{n,r}{\PP ^N_{j(\kappa )+1}}
\end{diagram}
$$
\end{figure}

and $N[K]$ is $\Sigma ^1_0$-correct for $\lambda$ in $V[G]$. (See
theorem~\ref{hauser}).
\end{itemize}

\medskip
\noindent
{\bf Proof of Theorem~\ref{presunf}:} Working in $N[K]$,
let $$({\frak A},\in ,\tilde S) =
\Bigl (V[K]_{j(\kappa )},\in ,j[(\dot S^M)^H]\Bigr ).$$
It is easy to check that $j[(\dot S^M)^H] =
j[(\dot S^M)^G]$. Then, suppose that
$$N[K]\models \Bigl [ ({\frak A},\in ,\tilde S)\models \exists x\varphi
(x,a_1,\dots ,a_n)\Bigr ],$$
where $a_1,\dots ,a_n\in [{\cal R}(\kappa )]^{V[K]}.$
Then, since $\kappa$ is the critical point of $j$,
$$M[H]\models \Bigl [ ({\cal R}(\kappa ),\in ,(\dot S)^G)\models
\exists x\varphi (x,a_1,\dots ,a_n)\Bigr ].$$
But since whenever $\alpha \leq \kappa$, $M[H]_\alpha = V[G]_\alpha$,
we have that
$$V[G]\models \Bigl [ ({\cal R}(\kappa ),\in ,(\dot S)^G)\models
\exists x\varphi (x,a_1,\dots ,a_n)\Bigr ].$$

\medskip
\noindent
But this implies that $\kappa$ is unfoldable in the extension by $\PP$:
we have just provided the embedding (in $V[G]$) required by our
characterisation of unfoldability in Theorem~\ref{theorem:embed}.
\hfill $\Box_{\ref{presunf}}$

\noindent
The much nicer version of the theorem that would be obtained by removing
the lifting conditions is unfortunately not true, as demonstrated by the
previously given counterexample. Still, we think that those lifting conditions
are natural enough as conditions for preservation of unfoldability.

\bigskip

\section{Destroying or adding GCH at an unfoldable.}

\noindent
Preservation of large cardinal properties while forcing new subsets of them
has played an important role in the development of iterated forcing notions.
We look next at iterations that preserve unfoldability of a cardinal
$\kappa$ while adding large numbers of subsets. Kai Hauser proved
in~[Ha 91] that $\Pi ^1_n$-indescribable cardinals
$\kappa$ are preserved
under iterations that add Cohen subsets of all inaccessible cardinals below
$\kappa$. We are interested in studying the consistency strength of
the failure of GCH at an unfoldable. Since we will be using iterations
of forcings that add new Cohen subsets at various places, we begin by
noticing the following useful fact.

\begin{fact}
If the $\lambda$-unfoldability of $\kappa$ is not destroyed by adding one Cohen
subset to $\kappa$, then it is not destroyed by adding
any number $\mu >\kappa$ of Cohen subsets to $\kappa$. \hfill $\Box$
\end{fact}

\noindent
We also remind the reader of the definition of ineffable cardinals.

\begin{definition}
A cardinal $\theta$ is {\bf ineffable} iff for any sequence
$\langle S_\alpha \subset \alpha | \alpha <\theta \rangle$ there exists
a `coherence' set $S\subset \theta$ such that $\{\alpha <\theta |
S\cap \alpha =S_\alpha \}$ is stationary in $\theta$.
\end{definition}

\noindent
The consistency strength of ineffables is way below that of Ramseys, but
much stronger than that of weakly compacts: if $\theta$ is ineffable, then
there are $\theta$ weakly compacts below it. On the other hand, ineffables
are still in the `small team' among large cardinals in the sense that they
do not imply the existence of $0^\sharp$: they actually relativise to $L$.

\noindent
Levinski's theorem (Theorem~\ref{levinski} here) directly provides us
with a way of {\bf adding} GCH at unfoldable cardinals.

\begin{theorem}\label{ineffgchunf}
Let $\kappa$ be ineffable. Then
the iteration of collapses that adds GCH defined in Theorem~\ref{levinski}
preserves the existence of unfoldables in ${\cal R}(\kappa )$.
\end{theorem}

\noindent
{\bf Proof:} by Theorem~\ref{levinski}, the ineffability of $\kappa$ is
preserved by $\PP _\kappa$, and GCH is added all the way beneath $\kappa$.
Then, by theorem~\ref{theorem:subtle}, there must be unfoldables in
the sense of ${\cal R}(\kappa )$. Of course, GCH holds at those
${\cal R}(\kappa )$-unfoldables.
\hfill $\Box$

\medskip
\noindent
We provide next an upper bound for the consistency strength of the failure
of GCH at an unfoldable cardinal.

\begin{theorem}\label{ineffunfold}
Con$\big[$ZFC $+\exists \theta$ $(\theta$ is
subtle)$\big]$ implies Con$\big[$ZFC $+\exists \kappa$ $(\kappa$ is
unfoldable $+\mbox { } 2^{\kappa}>\kappa )\big]$.
\end{theorem}

\noindent
In some sense, this is not a very tight upper bound for the consistency
strength of the failure of
GCH at an unfoldable cardinal. But this is sharp enough to prove that
the failure of GCH at an unfoldable cardinal does not imply the existence of
$0^\sharp$. The sharper upper bound Con$\big[$ZFC
$+\exists \theta$ $(\theta$ is unfoldable)$\big]$ seems plausible, but would
be obtained in a totally different way. Yet the
proof using subtlety provides some additional insight into the relation
between unfoldables and subtles (and ineffables). On the other hand, the
next lemma is interesting independently of its connection to unfoldable
cardinals.

\medskip
\noindent
{\bf Proof of Theorem~\ref{ineffunfold}:} Start with a model of GCH that
has a subtle cardinal
$\theta$. Let $\PP _\theta$ be an iteration of length $\theta$, such that
the active stages correspond to the strongly inaccessible cardinals
$\alpha <\theta$,
at which $\alpha ^{++}$ new Cohen subsets are added. 
Formally, define $\PP _\theta$ as the following iteration: $\PP _0 =$
the trivial forcing; for limit $\alpha < \theta$, $\PP _\alpha$
is the direct limit
of $\langle \PP _\beta | \beta < \alpha \rangle$ if $\alpha$ is inaccessible,
the inverse limit otherwise. Working on $\PP _\alpha$, $\dot \QQ _\alpha \in
V^{\PP _\alpha}$ is such that $\forces _{\PP _\alpha}`\dot \QQ _\alpha =
\{ 0\}$ if $\alpha$ is not inaccessible, $\dot \QQ _\alpha = \dot
{Add}(\alpha ^{++},\alpha)$ otherwise'; finally, $\PP _{\alpha +1} = \PP
_\alpha * \dot \QQ _\alpha$.

\noindent
Our aim is to prove that the subtlety of $\theta$ is preserved by this
iteration: If this is the case, and $G$ is $\PP _\theta$-generic, then,
by Theorem~\ref{theorem:subtle}, we know that the model
$[{\cal R}(\theta )]^{V[G]}$ `thinks' that there are unfoldable cardinals.
But in particular, GCH has been destroyed at those unfoldables, and so
$[{\cal R}(\theta )]^{V[G]}\models `$ZFC $+\exists \kappa$ $(\kappa$ is
unfoldable $+\mbox { } 2^{\kappa}>\kappa )$'.

\begin{lemma}\label{subtlpres}
The subtlety of $\theta$ is preserved under the iteration of forcings
that add new Cohen subsets to inaccessibles {\bf below} $\theta$.
\end{lemma}

\noindent
{\bf Proof of the Lemma:} We need to prove that, in $V[G]$,
for every sequence $\langle S_\alpha |
\alpha <\kappa \rangle$ such that $S_\alpha \subset \alpha$, for each
$\alpha <\kappa$, and for any $C$ club in $\kappa$,
there are at least two elements $\beta <\gamma$ of $C$ such that
$S_\beta = S_\gamma \cap \beta$.

\noindent
First, observe that clubs in $\kappa$
have to be in the ground model: no unbounded subset has been added to
$\kappa$. So, fix $C$ a club in $\kappa$ (in $V$).

\noindent
On the other hand, the sequence $\langle S_\alpha | \alpha <\kappa \rangle$
consists (in general) of elements of the extension; choose canonical names
$\sigma _\alpha$ for each $S_\alpha$. This can be done in the following
more canonical way: Let
$$1\forces _{\PP _\theta} (\dot F \mbox { is a function }\wedge
\dot F(\alpha )\subset
\alpha ,\forall \alpha \in C).$$ We need to prove that the values of $\dot F$
cohere on at least two elements of $C$. Now, $\dot F(\alpha )$ depends only
on $\alpha$ many generics (added at each level up to and including
$\alpha$). For each $\alpha \in C$, choose a
$\PP _{\alpha +1}$-name $\tau _\alpha$, and choose a condition $p_\alpha \in
\PP _\theta$ such that $$p_\alpha \forces _{\PP _\theta}(\tau
_\alpha = \dot F(\alpha )).$$
The $\theta$-cc of $\PP _\theta$ implies that there must exist $\alpha <
\beta$ from $C$ such that $p_\alpha$ and $p_\beta$ are compatible. Without
loss of generality, these two conditions can be chosen in such a way that
the supremum of the support of $p_\alpha$ is $<\beta$. Moreover, without
loss of generality, and by the subtlety of $\theta$ in the ground
model, the restrictions $p_\alpha \restriction \alpha$ and
$p_\beta \restriction \beta$ can be made coherent, so that
$p_\alpha \restriction \alpha = p_\beta \restriction \alpha$.

\noindent
Thus equipped, it only remains to make the $\tau$'s coherent at two places
in $C$ as well. We cannot apply subtlety of $\theta$ in $V$ yet, since
the $\tau _\alpha$'s do not live on the $\alpha$'s. Nevertheless, as each
$\tau _\alpha$ only depends on $\alpha$ many of the generic subsets
of $\alpha$, we can apply subtlety to subsets $\tilde \tau _\alpha$
of $\alpha$. These subsets represent isomorphism types of the possible
$\tau _\alpha$'s.

\noindent
The two previous applications of subtlety (to make two $p_\alpha$'s
cohere and then to make two $\tilde \tau _\alpha$'s cohere) can really be made
simultaneously, by coding the corresponding pairs. So, we get that
$$\exists \alpha <\beta, \alpha ,\beta \in C (p_\beta \forces
_{\PP _\theta} \tilde \tau _\beta \cap \alpha = \tau _\alpha).$$
But then, $\theta$ is also subtle in the generic extension.
\hfill $\Box _{\ref{subtlpres}}$

\noindent
With this, we finish the proof that GCH is destroyed at unfoldables,
by looking at the model ${\cal R}(\theta )$.
\hfill $\Box _{\ref{ineffunfold}}$

\medskip
\noindent
{\bf Remark:} a crucial feature of the proof of Lemma~\ref{subtlpres} was the
boundedness of the forcing; the fact that nothing was added at stage
$\theta$. The lemma shows that the subtlety of $\theta$ plus the failure of GCH
at all inaccessibles {\bf below} $\theta$ is not stronger
consistencywise than the subtlety of $\theta$. The related result that
the subtlety of $\theta$ plus the failure of GCH {\bf at} $\theta$
is not stronger consistencywise than the subtlety of $\theta$ stems from the
fact that $\Pi ^1_1$ properties at $\theta$ are not destroyed by
$(<\theta )$-closed forcings. (This last fact is due to Silver and Kunen).

\bigskip
\noindent
{\bf Remark: the failure of GCH at a fixed unfoldable:} the previous proof
requires the existence of a subtle $\theta$ and provides unfoldable
cardinals beneath $\theta$ at which the GCH fails. A natural question is
whether the failure of GCH may be obtained at $\kappa$, starting {\bf from
hypotheses about $\kappa$,} and not about larger cardinals.

\medskip
\noindent
{\bf Conjecture} Con(ZFC + $\exists \kappa (\kappa$ strongly unfoldable))
$\Rightarrow$ Con(ZFC + $\exists \kappa (2^\kappa > \kappa \wedge \kappa$
unfoldable)).

\noindent
This would mean that the consistency strength of the failure of GCH at
an unfoldable $\kappa$ is not greater than the strong unfoldability
of $\kappa$. This would sharpen the previous theorem.

\noindent
Here are some reasons for the previous conjecture: To start,
fix $V=L$, and $\kappa$ a $\lambda$-strongly unfoldable
cardinal. We plan to destroy GCH at $\kappa$ while proving that the
$\lambda$-unfoldability of $\kappa$ is preserved in the corresponding
forcing extension.
The natural forcing would be an iteration of length $\kappa +1$,
where the active
stages correspond to the strongly inaccessible cardinals $\alpha <\kappa$,
at which $\alpha ^{++}$ new Cohen subsets are added. Then we would use the
observation in the last section about the existence of a lifting for the
corresponding embeddings
to prove that the $\lambda$-unfoldability of $\kappa$ is preserved.

\medskip
\noindent
Formally, define $\PP _{\kappa +1}$ as the following iteration: $\PP _0 =$
the trivial forcing; for limit $\alpha$, $\PP _\alpha$ is the direct limit
of $\langle \PP _\beta | \beta < \alpha \rangle$ if $\alpha$ is inaccessible,
the inverse limit otherwise. Working on $\PP _\alpha$, $\dot \QQ _\alpha \in
V^{\PP _\alpha}$ is such that $\forces _{\PP _\alpha}`\dot \QQ _\alpha =
\{ 0\}$ if $\alpha$ is not inaccessible, $\dot \QQ _\alpha = \dot
{Add}(\alpha ^{++},\alpha)$ otherwise'; finally, $\PP _{\alpha +1} = \PP
_\alpha * \dot \QQ _\alpha$.

\noindent
Clearly, if $\mu \leq \kappa$ is inaccessible, then for all $\alpha <\mu$,
$|\PP _\alpha |<\mu$. So, if $\mu \leq \kappa$ is Mahlo, the forcing $\PP _\mu$
is $\mu$-cc, since then $\{\alpha <\mu | \PP _\alpha = \overrightarrow\lim
_{\eta <\alpha} \PP _\eta \}$ is stationary in $\mu$. Also, $\forces _{\PP
_\mu}`\mu$ is inaccessible', and for all $\alpha <\mu$, $\forces _{\PP
_\alpha}`\PP _{\alpha ,\mu}$ is $(<\nu)$-closed', where $\nu$ is the least
inaccessible $\geq \alpha$ and $\PP _{\alpha ,\mu}$ is the `remainder forcing'.

\medskip
\noindent
Fix $M$ a model of $ZF^-$ of size $\kappa$, where $\kappa \in M$ and
$M^{<\kappa}\subset M$, and get $N$ and $j:M\to N$ a $\lambda$-unfoldability
embedding as in Theorem \ref{theorem:embed}..

\medskip
\noindent
{\bf The obstruction:} To reach the conclusion that
$\forces _{\PP _{\kappa +1}}`\kappa$ is $\lambda$-unfoldable',
we would need to produce $G$ $\PP ^M_{\kappa +1}$-generic over $M$. A
natural way to obtain $G$ would be by using the `lifting conditions' mentioned
in the last section. These, nevertheless, require an amount of closure
for the portion of the forcing beyond $\kappa$ that is not yet met.

\bigskip

\appendix
\section{Class Forcing and eees.}

\medskip
\noindent
In this section, we consider what happens to models of $ZFC$ and their
extendibility properties after various sorts of Class Forcing are applied.
We first mention the observation made in [Vi $\infty$]

\begin{fact}
If ${\cal E}_M \not= 0$, and G is SET-generic over M, then, ${\cal E}_{M[G]}
\not= 0$.
\end{fact}

\noindent
When applying CLASS forcing, this situation may change dramatically: the fact
that $M[G]\models ZFC$ is no longer guaranteed. Also, the proof that
${\cal E}_{M[G]}\not= 0$ when ${\cal E}_M \not= 0$ uses in a nontrivial way
the fact that $G$ is $\PP$--generic, for some $\PP$ element of $M$.

\bigskip
\noindent
We now restrict ourselves to the countable case, and provide an application of
a theorem due to Aaron Beller from the famous book `Coding the Universe'
by Beller, Jensen and Welch~[BeJeWe 83]. Using Keisler and Morley's
[KeMo 68], there is no problem here with the existence of eees.
Nevertheless, the following fact gives insight into the situation
we get into when applying class forcing.

\begin{proposition} There exist countable $M\models ZFC +$ `unfoldable
cardinals exist', and $a\subset \omega$ such that $M[a]\models ZFC +$ `no
unfoldable cardinals exist'.
\end{proposition}

\noindent
{\bf Proof:} Remember that $M[a]$ denotes the minimum model that contains
$M\cup\{ a\}$. In [BeJeWe 83], Aaron Beller proves that given any $\beta
<\omega_1$, there is a real $a\subset \omega$ such that given any $\alpha
<\beta$, $L_\alpha [a]\not\models ZFC$. So, pick $\beta <\omega _1$ and
$a\subset \omega$ such that
$L_\beta$ is a limit of a chain of eees,
$L_\beta [a]\models ZFC$, and for any $\alpha
<\beta$, $L_\alpha [a]\not\models ZFC$. Clearly, $L_\beta [a]\models$
`no unfoldable cardinals exist.' But since $\beta$ was chosen as a limit of
a chain of eees, $L_\beta \models ZFC +$ `unfoldable cardinals exist', yet
the (class generic) extension $L_\beta [a]$ doesn't. \hfill $\Box$

\medskip
\noindent
We finish by posing the following

\begin{question} Is there a model of $ZFC+\mbox {`unfoldable
cardinals exist'}$, such that for some $a\subset \omega$, $M[a]\models
ZFC+\mbox {`no unfoldable cardinals exist'}$, yet $M[a]\models
ZFC+\mbox {`weakly compact cardinals exist'}$?
\end{question}

\bigskip
\noindent
{\Large \bf References.}
\begin{description}

\item[(BeJeWe 83)] {\capit Beller, A., Jensen, Welch}{ \it Coding the
Universe}, Cambridge Lecture Series, 1983.
\item[(Bo 74)] {\capit Boos, W.}{ \it Boolean Extensions which efface
the Mahlo property}, The Journal of Symbolic Logic, vol. 39, no 2, June 74,
pages 254-268.
\item[(En 84)] {\capit Enayat, A.}{ \it On certain Elementary Extensions of
Models of Set Theory}, Trans. Amer. Math. Soc., 1984.
\item[(En $\infty$)] {\capit Enayat, A.}{ \it Counting Countable
Models of Set Theory}, preprint.
\item[(Ha 91)] {\capit Hauser, K.}{ \it Indescribable Cardinals and
Elementary Embeddings}, The Journal of Symbolic Logic, vol. 56, no 2, June 91,
pages 439-457.
\item[(Je 74)] {\capit Jensen, R.}{ \it Measurable Cardinals and the GCH},
in T. Jech, ed., Axiomatic Set Theory, Proc. Symp. Pure Math., Vol. 13,
II (American Mathematical Society, Providence, RI, 1974), 175-178
\item[(Ka 83)] {\capit Kaufmann, M.}{ \it Blunt and Topless Extensions of
Models of Set Theory}, J. Symb. Logic, 48, 1983, 1053-1073.
\item[(KeMo 68)] {\capit Keisler, H.J., Morley, M.}{ \it Elementary
Extensions of Models
of Set Theory.} Israel J. of Math., vol. 6, 1968.
\item[(KeSi 70)] {\capit Keisler, H.J., Silver, J.} { \it
End Extensions of Models of Set Theory.}, Proc. Symp. Pure Math. 13
1970 (177-187).
\item[(La 78)] {\capit Laver, R.} { \it
Making the Supercompactness of $\kappa$ indestructible
under $\kappa$-directed closed forcing}, Israel Journal of Mathematics, 29,
1978 (385-388).
\item[(Le 95)] {\capit Levinski, J.P.}{ \it Filters and Large Cardinals},
Annals of Pure and Applied Logic, vol. 72, 1995, 177-212.
\item[(Vi $\infty$)] {\capit Villaveces, A.}{ \it Heights of Models of $ZFC$
and the Existence of End Elementary Extensions}, preprint.
\end{description}
\end{document}